\documentstyle[amscd,amssymb,verbatim,epsf]{amsart}

\begin{document}
\theoremstyle{plain}
\newtheorem{Thm}{Theorem}
\newtheorem{Cor}{Corollary}
\newtheorem{Ex}{Example}
\newtheorem{Con}{Conjecture}
\newtheorem{Main}{Main Theorem}
\newtheorem{Lem}{Lemma}
\newtheorem{Prop}{Proposition}

\theoremstyle{definition}
\newtheorem{Def}{Definition}
\newtheorem{Note}{Note}

\theoremstyle{remark}
\newtheorem{notation}{Notation}
\renewcommand{\thenotation}{}

\errorcontextlines=0
\numberwithin{equation}{section}
\renewcommand{\rm}{\normalshape}%

\title[K\"{a}hler Metric]%
   {An Indefinite K\"{a}hler Metric on the Space of Oriented Lines}
\author{Brendan Guilfoyle}
\address{Brendan Guilfoyle\\
          Department of Mathematics and Computing \\
          Institute of Technology, Tralee \\
          Clash \\
          Tralee  \\
          Co. Kerry \\
          Ireland.}
\email{brendan.guilfoyle@@ittralee.ie}
\author{Wilhelm Klingenberg}
\address{Wilhelm Klingenberg\\
 Department of Mathematical Sciences\\
 University of Durham\\
 Durham DH1 3LE\\
 United Kingdom.}
\email{wilhelm.klingenberg@@durham.ac.uk }

\date{March 23rd, 2004}

\begin{abstract}
The total space of the tangent bundle of a K\"ahler manifold admits a
canonical K\"ahler structure. Parallel translation identifies the
space ${\Bbb{T}}$ of oriented affine lines in ${\Bbb{R}}^3$ with
the tangent bundle of $S^2$. Thus, the round metric on $S^2$ induces a
K\"ahler structure on ${\Bbb{T}}$ which turns out to have a metric of
neutral signature.  It is shown that the
isometry group of this metric is isomorphic to 
the isometry group of the Euclidean metric on ${\Bbb{R}}^3$. 

The geodesics of this metric are either planes or helicoids in
${\Bbb{R}}^3$. The signature of the metric induced on a surface $\Sigma$ in
${\Bbb{T}}$ is determined by the degree of twisting of the associated
line congruence in ${\Bbb{R}}^3$, and we show that, for $\Sigma$
Lagrangian, the metric is either Lorentz or totally null. For such
surfaces it is proven that the Keller-Maslov index counts the number
of isolated complex points of ${\Bbb{J}}$ inside a closed curve on $\Sigma$.

\end{abstract}

\maketitle

\section{Introduction}
The total space of the tangent bundle of a K\"ahler $2n$-manifold ($N$,$j$,$\omega$,$g$),
denoted $TN$, itself carries a canonical K\"ahler structure
(${\Bbb{J}}$,$\Omega$,${\Bbb{G}}$). If $n=1$, this metric
${\Bbb{G}}$ has neutral signature ($++--$) and is conformally flat if
and only if $g$ is of constant curvature. 

The minitwistor correspondence identifies the tangent bundle of the
2-sphere with the space ${\Bbb{T}}$ of oriented affine lines in
${\Bbb{R}}^3$. Thus, the round metric on $S^2$ induces a K\"ahler
structure on ${\Bbb{T}}$ with the above properties. Since the
metric ${\Bbb{G}}$ has not been studied before in this context, the
purpose of this paper is to consider its significance for
differential geometry in ${\Bbb{R}}^3$.

The complex structure on ${\Bbb{T}}$
has proved crucial in understanding static monopoles which arise in
theoretical physics \cite{hitch}. Indeed, as early as 1866,
Weierstrass \cite{weier} used this complex structure to construct
minimal surfaces in ${\Bbb{R}}^3$, while Whittaker used it to find
solutions of the Laplace equation \cite{whitt}.  The symplectic
structure $\Omega$ is equivalent to the canonical symplectic structure
on $T^*S^2$, which comes to the fore in geometric optics (see, for
example, Arnold \cite{arn}).  

The interpretation we give to ${\Bbb{G}}$ is that of angular momentum
of the Jacobi fields along an oriented line in ${\Bbb{R}}^3$. This is
invariant under the group of rotations and translations, and
so the isometry group of ${\Bbb{G}}$ contains the
Euclidean group. In fact, we show that:

\vspace{0.1in}
\noindent{\bf Theorem 1.}
{\it
The  isometry group of the metric
${\Bbb{G}}$ on ${\Bbb{T}}$ is isomorphic to the
Euclidean isometry group.  
}
\vspace{0.1in}

\noindent This is analogous to the classical result on the automorphism
group of the Study sphere \cite{study}.

A curve in ${\Bbb{T}}$ generates a ruled surface in ${\Bbb{R}}^3$. We
prove that the geodesics of the metric ${\Bbb{G}}$ give particularly
simple ruled surfaces:
\vspace{0.1in}

\noindent{\bf Theorem 2.}
{\it
A geodesic of ${\Bbb{G}}$ is either a plane or a helicoid in
${\Bbb{R}}^3$, the former in the case when the geodesic is null, the
latter when it is non-null. 
}
\vspace{0.1in}

Of particular interest to us are surfaces $\Sigma\subset{\Bbb{T}}$,
classically referred to as {\it line congruences}. By construction,
$\Sigma$ is {\it Lagrangian} iff the associated line congruence in
${\Bbb{R}}^3$ has zero twist, that is, it is integrable in the sense 
of Frobenius and there exists orthogonal surfaces to the lines in ${\Bbb{R}}^3$. On the other hand, $\Sigma$ is {\it complex} iff the shear of
the line congruence vanishes.

We prove that the signature of the metric induced on
$\Sigma$ by ${\Bbb{G}}$ is determined by the ratio of twist to shear
of the congruence: 

\vspace{0.1in}
\noindent{\bf Theorem 3.}
{\it
The induced metric on a surface $\Sigma$ in ${\Bbb{T}}$  is 
Lorentzian (totally null, Riemannian) if and only if
$\lambda^2-|\sigma|^2<0(=0,>0)$,
where $\lambda$ and $\sigma$ are the twist and shear of the line
congruence $\Sigma$.
}
\vspace{0.1in}

As a consequence  the induced metric on a Lagrangian
surface is either Lorentz or totally null. The null directions on the 
orthogonal surface $S$ in ${\Bbb{R}}^3$ are the eigen-directions of
the second fundamental form, while the totally null points are
umbilical points on $S$.  

The Keller-Maslov index \cite{arn} associates an integer to 
closed oriented curves on a Lagrangian surface $\Sigma$ in a
symplectic 4-manifold. This is just the degree of
the Gauss mapping which takes a point on the curve to the tangent
plane of $\Sigma$, considered as a point in the Lagrangian
Grassmanian $\Lambda=U(2)/O(2)=S^1$. 

The reduction of $\mbox{GL}(2,{\Bbb{C}})$ to $U(2)$ is achieved by
choosing a complex structure which is tamed by $\Omega$, that is
${\Bbb{G}}$ is positive definite. In our case, ${\Bbb{J}}$ is not
tamed by the symplectic structure and we look at the metric reduction
of $\mbox{GL}(2,{\Bbb{C}})$ to $U(1,1)$. Since $U(1,1)/O(1,1)=S^1$, we
can still define an index on Lagrangian surfaces and we prove that:

\vspace{0.1in}
\noindent{\bf Theorem 4.}
{\it
The Keller-Maslov index of a curve on a Lagrangian surface in ${\Bbb{T}}$
counts the number of isolated complex points of ${\Bbb{J}}$ inside the curve.
}
\vspace{0.1in}

This paper is laid out as follows. The next section contains the
general definition of the K\"ahler structures on tangent bundles. In
section 3 we describe the space ${\Bbb{T}}$ of oriented affine lines in
${\Bbb{R}}^3$ and give the geometric interpretation of the neutral metric on ${\Bbb{T}}$. 

Section 4 contains the proof that the isometry group of
${\Bbb{G}}$ is isomorphic to the Euclidean group, while section 5
investigates the geodesics in ${\Bbb{T}}$. In section 6 we turn to the
2-dimensional submanifolds of ${\Bbb{T}}$ and their induced
geometry. In the final section we relate the Keller-Maslov index on
Lagrangian surfaces to the index of isolated complex points.

\section{Canonical K\"ahler Metric on Tangent Bundles}

Throughout this paper we work solely with smooth maps and manifolds. Let $M$ be a smooth even-dimensional manifold. An {\it almost complex
structure} on $M$ is an endomorphism ${\Bbb{J}}:T_pM\rightarrow T_pM$
at each $p\in M$ satisfying ${\Bbb{J}}\circ{\Bbb{J}}=
-\mbox{Id}$. If, in addition, ${\Bbb{J}}$ satisfies an
integrability condition, it is said to be a {\it complex
structure}. 

A {\it symplectic structure} on $M$ is a closed non-degenerate 2-form
$\Omega$. It is {\it compatible} with the almost complex structure if
$\Omega({\Bbb{J}}\cdot,{\Bbb{J}}\cdot)=\Omega(\cdot,\cdot)$.

Given a complex structure and compatible symplectic structure 
define a symmetric non-degenerate 2-tensor ${\Bbb{G}}:T_pM\times
T_pM\rightarrow {\Bbb{R}}$  by
${\Bbb{G}}(\cdot,\cdot)=\Omega({\Bbb{J}}\cdot,\cdot)$.
The manifold $M$ with this triple of structures is called a
{\it K\"ahler} manifold. While ${\Bbb{G}}$ is a metric, we make no
assumption that it is positive definite. 

Given a K\"ahler $2n$-manifold ($N$,$j$,$\omega$,$g$) we construct
a canonical K\"ahler structure (${\Bbb{J}}$,$\Omega$,${\Bbb{G}}$) on the tangent bundle $TN$ as follows. The
Levi-Civita connection associated with $g$ splits the tangent bundle
$TTN\equiv TN\oplus TN$ and the complex structure is
defined to be ${\Bbb{J}}=j\oplus j$. The integrability of
${\Bbb{J}}$ follows from the integrability of $j$ \cite{lich}.
To define the symplectic form, consider the metric $g$ as a mapping
from $TN$ to $T^*N$ and pull back the canonical symplectic form
$\Omega^*$ on $T^*N$ to $\Omega$ on $TN$. Finally, the metric is
defined as above by
${\Bbb{G}}(\cdot,\cdot)=\Omega({\Bbb{J}}\cdot,\cdot)$. The
triple (${\Bbb{J}}$, $\Omega$, ${\Bbb{G}}$) determine a K\"ahler
structure on $TN$. 

We now find local coordinate descriptions of the above. Let $\xi^k$ be
local holomorphic coordinates on an open set of $N$ such that the
metric has the real potential $Q$:
\[
g=\partial_k\bar{\partial}_lQ\; d\xi^k\otimes d\bar{\xi}^l
\qquad\qquad \mbox{where}\qquad \partial_i=\frac{\partial}{\partial \xi^i}.
\] 
Introduce coordinates on $T^*N$ by identifying
$(\xi^k,\alpha_l)$ with $\alpha_kd\xi^k+\bar{\alpha}_kd\bar{\xi}^k$, 
and coordinates on $TN$ by identifying$(\xi^k,\eta^l)$
with $\eta^k\partial_k+\bar{\eta}^k\bar{\partial}_k$. 

The canonical symplectic
structure on $T^*N$ is $\Omega^*=d\alpha_k\wedge
d\xi^k+d\bar{\alpha}_k\wedge d\bar{\xi}^k$.
Considered as a map between $TN$ and $T^*N$, the metric $g$ takes
($\xi^k$,$\eta^l$) to
($\xi^k$,$\alpha_l=\bar{\eta}^m\bar{\partial}_m\partial_lQ
$). Thus, pulling back the canonical symplectic form we get
\[
\Omega={\Bbb{R}}\mbox{e}\left(\partial_k\bar{\partial}_lQ \;d\eta^k\wedge d\bar{\xi}^l
   +\eta^m\partial_m\partial_k\bar{\partial}_lQ\;
  d\xi^k\wedge d\bar{\xi}^l\right).
\]

\begin{Prop}
The eigenspaces of ${\Bbb{J}}$ are:
\begin{equation}\label{e:gencomp}
{\Bbb{J}}\left(\frac{\partial}{\partial
    \xi^k}\right)=i\frac{\partial}{\partial \xi^k} \qquad \qquad
  {\Bbb{J}}\left(\frac{\partial}{\partial
    \eta^k}\right)=i\frac{\partial}{\partial \eta^k}.
\end{equation}
The complex structure ${\Bbb{J}}$ and the symplectic structure
$\Omega$ are compatible.
\end{Prop}

\begin{pf}
This follows from the fact that the metric connection of $g$ is
compatible with the complex structure $j$. 
\end{pf}

We now specialise to the case $n=1$, where more can be said:

\begin{Prop}
Let ($N$, $g$) be a riemannian 2-manifold. Then ($TN$,${\Bbb{J}}$, $\Omega$,
${\Bbb{G}}$) is a K\"ahler manifold. The metric ${\Bbb{G}}$ has
neutral signature ($++--$) and is scalar-flat. Moreover, ${\Bbb{G}}$ is K\"ahler-Einstein iff $g$ is flat, and
${\Bbb{G}}$ is conformally flat iff $g$ is of constant curvature.
\end{Prop}

\begin{pf}
Choose isothermal coordinates $\xi$ on $N$ so that
$ds^2=e^{2u}d\xi d\bar{\xi}$, and corresponding coordinates
($\xi$,$\eta$) on $TN$, as above. In such a coordinate system, the symplectic 2-form is

\begin{equation}\label{e:riemsymp}
\Omega=2{\Bbb{R}}\mbox{e}\left(e^{2u}d\eta\wedge d\bar{\xi}+\eta\partial(e^{2u})d\xi\wedge d\bar{\xi}\right),
\end{equation}
and the metric ${\Bbb{G}}$ is 
\begin{equation}\label{e:riemmetric}
{\Bbb{G}}=2{\Bbb{I}}\mbox{m}\left(e^{2u}d\eta d\bar{\xi}+\eta\partial(e^{2u})d\xi d\bar{\xi}\right).
\end{equation}
A direct computation shows that ${\Bbb{G}}$ is neutral and
scalar-flat.

The only non-vanishing component of the Ricci tensor is
$R_{\xi\bar{\xi}}=-4\partial\bar{\partial}u.$ Since the Gauss
curvature of $N$ is $\kappa=-4e^{-2u}\partial\bar{\partial}u$,
${\Bbb{G}}$ is K\"ahler-Einstein iff $\kappa=0$. 

Furthermore, the only non-vanishing component of the conformal
curvature tensor is
\[
C_{\xi\bar{\xi}\xi\eta}=\frac{1}{2}e^{2u}(\eta\partial \kappa+\bar{\eta}\bar{\partial} \kappa),
\]
so ${\Bbb{G}}$ is conformally flat iff $\kappa$ is constant.
\end{pf}

\section{The K\"{a}hler Metric on ${\Bbb{T}}$}

Consider the space ${\Bbb{T}}$ of oriented (affine) lines in
${\Bbb{R}}^3$. Each oriented line can be uniquely described by two vectors: the
unit direction vector of the line and the perpendicular distance
vector of the line from the origin. If we think of the direction
vector as a point $\xi$ on the 2-sphere, the perpendicular distance
vector gives a tangent vector to $S^2$ at $\xi$.
Thus ${\Bbb{T}}$ is diffeomorphic to the tangent bundle to $S^2$, giving it the structure of a smooth 4-manifold \cite{hitch}.

Take local complex coordinates $\xi$  on  $S^2$ by
stereographic projection from the South pole onto the plane through
the equator. These coordinates yield canonical coordinates
($\xi$,$\eta$) on ${\Bbb{T}}= TS^2$, as in the last section.

\begin{Prop}\cite{gak1}
Consider the map $\Phi:{\Bbb{T}}\times{\Bbb{R}}
\rightarrow {\Bbb{R}}^3$ defined by   
\begin{equation}\label{e:coord1}
\Phi_z=\frac{2(\eta-\bar{\eta}\xi^2)+2\xi(1+\xi\bar{\xi})r}{(1+\xi\bar{\xi})^2}
\end{equation}
\begin{equation}\label{e:coord2}
\Phi_t=\frac{-2(\eta\bar{\xi}+\bar{\eta}\xi)+(1-\xi^2\bar{\xi}^2)r}{(1+\xi\bar{\xi})^2},
\end{equation}
Then $\Phi$ takes a point $(\xi,\eta)\in T_\xi S^2$ and $r\in{\Bbb{R}}$, to a point $(\Phi_z,\Phi_t)\in{\Bbb{C}}\oplus{\Bbb{R}}={\Bbb{R}}^3$
which is a distance $r$ along the line with direction $\xi$ and
minimal distance from the origin given by $(\xi,\eta)$.
\end{Prop} 

A null frame at a point in ${\Bbb{R}}^3$ is a trio of vectors $e_{(0)}$,
$e_{(+)}$, $e_{(-)}\in{\Bbb{C}}\otimes T{\Bbb{R}}^3$ such that
$e_{(0)}=\overline{e_{(0)}}$,$e_{(+)}=\overline{e_{(-)}}$, $e_{(0)}\cdot
e_{(0)}=e_{(+)}\cdot e_{(-)}=1$, and $e_{(0)}\cdot e_{(+)}=0$, where the
Euclidean inner product is extended bilinearly over ${\Bbb{C}}$.

The derivative of the map $\Phi$ has the following description:

\begin{Prop}
The derivative
$D\Phi:T_{(\xi,\eta,r)}{\Bbb{T}}\times{\Bbb{R}}\rightarrow T_{\Phi(\xi,\eta,r)}{\Bbb{R}}^3$ is given by 
\begin{equation}\label{e:dphi1}
D\Phi_{(\xi,\eta,r)}\left(\frac{\partial}{\partial
    \xi}\right)=\left(r-\frac{2\bar{\xi}\eta}{1+\xi\bar{\xi}}\right)\frac{\sqrt{2}}{1+\xi\bar{\xi}}e_{(+)}-\frac{2\bar{\eta}}{(1+\xi\bar{\xi})^2}e_{(0)}
\end{equation}
\begin{equation}\label{e:dphi2}
D\Phi_{(\xi,\eta,r)}\left(\frac{\partial}{\partial
    \eta}\right)=\frac{\sqrt{2}}{1+\xi\bar{\xi}}e_{(+)}
\qquad\qquad
D\Phi_{(\xi,\eta,r)}\left(\frac{\partial}{\partial
    r}\right)=e_{(0)},
\end{equation}
where we have introduced the null frame
\[
e_{(0)}= \frac{2\xi}{1+\xi\bar{\xi}}\frac{\partial}{\partial z}
     +\frac{2\bar{\xi}}{1+\xi\bar{\xi}}\frac{\partial}{\partial \bar{z}}
  + \frac{1-\xi\bar{\xi}}{1+\xi\bar{\xi}}\frac{\partial}{\partial t}
\]
\[
e_{(+)}= \frac{\sqrt{2}}{1+\xi\bar{\xi}}\frac{\partial}{\partial z}
     -\frac{\sqrt{2}\;\bar{\xi}^2}{1+\xi\bar{\xi}}\frac{\partial}{\partial \bar{z}}
  - \frac{\sqrt{2}\;\bar{\xi}}{1+\xi\bar{\xi}}\frac{\partial}{\partial t}.
\]
Here and throughout we use $z=x^1+ix^2$ and $t=x^3$, where $x^1$, $x^2$
and $x^3$ are Euclidean coordinates on ${\Bbb{R}}^3$. The unit vector
$e_{(0)}$ determines the direction of the line given by ($\xi$, $\eta$).
\end{Prop}

The map $D\Phi$ gives the identification of tangent
vectors to ${\Bbb{T}}$ at a line $\gamma$ with the Jacobi fields
orthogonal to the line in ${\Bbb{R}}^3$. A Jacobi field along a line
$\gamma$ in ${\Bbb{R}}^3$ is a vector field $X$ along the line which
satisfies the equation
\[
\nabla_{\dot{\gamma}}\nabla_{\dot{\gamma}}X=0
\]
Choosing an affine parameter $r$ along the line, this has solution
$X=X_1+rX_2$, for constant vector fields $X_1$ and $X_2$ along $\gamma$. The
Jacobi fields that are orthogonal to $\gamma$ form a 4-dimensional
vector space, which $D\Phi$ identifies with the tangent space
$T_\gamma{\Bbb{T}}$.

From the maps (\ref{e:dphi1}) and
(\ref{e:dphi2}), along with equation (\ref{e:gencomp}), the complex
structure $\Bbb{J}$ acts on the Jacobi fields by multiplication of
$e_{(+)}$ by $i$. This is equivalent to rotation about $e_{(0)}$
through 90$^0$. This yields (see \cite{hitch}):

\begin{Prop}
The complex structure $\Bbb{J}$ on ${\Bbb{T}}$ defined in the last
section is given by rotation of the Jacobi fields through 90$^o$ about
the corresponding line in ${\Bbb{R}}^3$. 
\end{Prop}

\begin{Prop} 
The symplectic 2-form $\Omega$ on ${\Bbb{T}}$,  determined by the
round metric on $S^2$, is given by 
\begin{equation}\label{d:symp}
\Omega_{(\xi,\eta)}({\Bbb{X}},{\Bbb{Y}})=\left<D\Phi({\Bbb{X}}),\nabla_{(0)}D\Phi({\Bbb{Y}})\right> -\left<D\Phi({\Bbb{Y}}),\nabla_{(0)}D\Phi({\Bbb{X}})\right>,
\end{equation}
where ${\Bbb{X}},{\Bbb{Y}}\in T_{(\xi,\eta)}{\Bbb{T}}$,
$<\cdot,\cdot>$ is the Euclidean metric on ${\Bbb{R}}^3$ and
$\nabla_{(0)}$ is the covariant derivative in the $\xi$ direction.
\end{Prop}
\begin{pf}
Taking (\ref{d:symp}) as the definition of $\Omega$, we find, for
example, using (\ref{e:dphi1}) and (\ref{e:dphi2}), that
\begin{align}
\Omega\left(\frac{\partial}{\partial\xi},\frac{\partial}{\partial\bar{\xi}}\right)
=&
g\left(\frac{\sqrt{2}}{1+\xi\bar{\xi}}\left(r-\frac{2\bar{\xi}\eta}{1+\xi\bar{\xi}}\right)e_{(+)},\frac{\sqrt{2}}{1+\xi\bar{\xi}}e_{(-)}\right)\nonumber\\
&\qquad-g\left(\frac{\sqrt{2}}{1+\xi\bar{\xi}}\left(r-\frac{2\xi\bar{\eta}}{1+\xi\bar{\xi}}\right)e_{(-)},\frac{\sqrt{2}}{1+\xi\bar{\xi}}e_{(+)}\right)\nonumber\\
=&\frac{2}{(1+\xi\bar{\xi})^2}\left(r-\frac{2\bar{\xi}\eta}{1+\xi\bar{\xi}}\right)-\frac{2}{(1+\xi\bar{\xi})^2}\left(r-\frac{2\xi\bar{\eta}}{1+\xi\bar{\xi}}\right)\nonumber\\
=&\frac{4(\xi\bar{\eta}-\bar{\xi}\eta)}{(1+\xi\bar{\xi})^3},\nonumber
\end{align}
and this agrees with the canonical symplectic structure given in
(\ref{e:riemsymp}) in the case when $g$ is the round metric on $S^2$
with $e^{2u}=2(1+\xi\bar{\xi})^{-2}$. Similarly for the other components of the symplectic form.
\end{pf}

The metric ${\Bbb{G}}$ on ${\Bbb{T}}$ is defined by
${\Bbb{G}}(\cdot,\cdot)\equiv\Omega(\cdot,{\Bbb{J}}\cdot)$. By
Proposition 2 we have that:

\begin{Prop}
The metric ${\Bbb{G}}$ on ${\Bbb{T}}$ is a neutral K\"ahler metric
which is conformally flat, with zero scalar curvature, but is not Einstein.
\end{Prop}

The local expression for the metric is
\begin{equation}
{\Bbb{G}}=\frac{2i}{(1+\xi\bar{\xi})^2}\left(
  d\eta d\bar{\xi}-d\bar{\eta} d\xi
   +\frac{2(\xi\bar{\eta}-\bar{\xi}\eta)}{1+\xi\bar{\xi}}d\xi d\bar{\xi}
\right),
\end{equation}
and it can be given the following geometric interpretation:

\begin{Prop}
The length of ${\Bbb{X}}\in T_\gamma{\Bbb{T}}$ with respect to ${\Bbb{G}}$ is the
{\it angular momentum} about $\gamma$ of the line determined by the
Jacobi field associated to ${\Bbb{X}}$. 
\end{Prop}
\begin{pf}
A direct computation using $D\Phi$ shows that the length of 
${\Bbb{X}}=X_1+rX_2$ is the oriented area of the parallelogram spanned
by $X_1$ and $X_2$, i.e. ${\Bbb{G}}({\Bbb{X}},{\Bbb{X}})=(X_1\times
X_2)\cdot e_{(0)}$. 
\end{pf}

The symplectic 2-form on ${\Bbb{T}}$ is globally exact $\Omega=d\Theta$. This global
1-form, locally pulled back in the above coordinates, is
\begin{equation}\label{e:contact}
\Theta=\frac{2\bar{\eta}d\xi}{(1+\xi\bar{\xi})^2}
    +\frac{2\eta d\bar{\xi}}{(1+\xi\bar{\xi})^2}.
\end{equation} 
We discuss this further in section 6. The K\"{a}hler potential of the metric is
\[
\Upsilon=\frac{2i(\xi\bar{\eta}
  -\bar{\xi}\eta)}{1+\xi\bar{\xi}}.
\]

\section{The Isometry Group of ${\Bbb{G}}$}

The group of fibre-preserving holomorphic automorphisms
of ${\Bbb{T}}$ is isomorphic to $\mbox{PSL}(2,{\Bbb{C}}) \ltimes {\Bbb{C}}^3$,
where the action on ${\Bbb{T}}$ is given by
\begin{equation}\label{e:trans}
\xi\rightarrow\xi'=\xi  \qquad 
\eta\rightarrow \eta'=\eta+a_1+b_1\xi-c_1\xi^2,
\end{equation}
and
\begin{equation}\label{e:rot}
\xi\rightarrow\xi'=\frac{a_2\xi+b_2}{c_2\xi+d_2}  \qquad 
\eta\frac{\partial}{\partial \xi}\rightarrow \eta'\frac{\partial}{\partial \xi'}=\frac{\eta}{(c_2\xi+d_2)^2}\frac{\partial}{\partial \xi'},
\end{equation}
for $a_1,a_2,b_1,b_2,c_1,c_2,d_2\in{\Bbb{C}}$ with $a_2d_2-b_2c_2=1$. 

Here, $A \in \mbox{PSL}(2,{\Bbb{C}})$ acts on the $2 \times 2$ 
complex symmetric matrices (to be identified with 
${\Bbb{C}}^3$) by $M \to A M A^{T}$. This is the same action as
that of  $\mbox{PSL}(2,{\Bbb{C}})$ on quadratic holomorphic transformations
of ${\Bbb{T}}$ as above. 

The identity component of the Euclidean isometry group is double
covered by $\mbox{SU}(2)\ltimes{\Bbb{R}}^3$, which is a real form
of the complex Lie group $\mbox{PSL}(2,{\Bbb{C}}) \ltimes
{\Bbb{C}}^3$, where the transformations (\ref{e:trans}) and
(\ref{e:rot}) are restricted to those with
$a_1=\overline{c_1}$, $b_1=\overline{b_1}$, $a_2=\overline{a_2}=d_2$
and $b_2=-\overline{c_2}$.

Since the above construction of ${\Bbb{G}}$ is invariant under
Euclidean motions, it is clear that the isometry group of ${\Bbb{G}}$
contains the Euclidean group of translations and rotations. We now
prove that ${\Bbb{G}}$ admits no other continuous isometries:

\begin{Thm}
The isometry group of the metric
${\Bbb{G}}$ on ${\Bbb{T}}$ is isomorphic to the
Euclidean isometry group.  
\end{Thm}
\begin{pf}
The Lie group of Euclidean motions is a subgroup of the
Lie group of isometries of ${\Bbb{G}}$. This can be proved in local
coordinates, or, by noting that the angular momentum of a Jacobi field
is invariant under Euclidean motions.

We prove in the following proposition that the associated Lie algebras
are isomorphic, and in particular, have the same dimension. Thus the
connected component of the identity of the isometry groups of 
${\Bbb{G}}$ is isomorphic to the identity component of the Euclidean
group.  In addition, both groups are double covered by $\mbox{SU}(2)\ltimes
{\Bbb{R}}^3$ and the result follows.
\end{pf}

\begin{Prop}
The Killing vectors of ${\Bbb{G}}$ form a 6-parameter Lie algebra given by
\[
{\Bbb{K}}={\Bbb{K}}^\xi\frac{\partial}{\partial \xi}
  +{\Bbb{K}}^\eta\frac{\partial}{\partial \eta}
  +{\Bbb{K}}^{\bar{\xi}}\frac{\partial}{\partial \bar{\xi}}
  +{\Bbb{K}}^{\bar{\eta}}\frac{\partial}{\partial \bar{\eta}},
\]
with
\[
{\Bbb{K}}^\xi=\alpha+2ai\xi+\bar{\alpha}\xi^2 \qquad
{\Bbb{K}}^\eta=2(ai+\bar{\alpha}\xi)\eta+\beta+b\xi-\bar{\beta}\xi^2,
\]
where $\alpha,\beta\in{\Bbb{C}}$ and $a,b\in{\Bbb{R}}$. The Killing
vectors given by $\alpha$ and $a$ generate infinitesimal rotations
while those given by $\beta$ and $b$ generate infinitesimal translations.

\end{Prop}
\begin{pf}
A vector field ${\Bbb{K}}^i$ is a Killing vector for ${\Bbb{G}}$ if and only if
\begin{equation}\label{e:killing}
{\Bbb{K}}^i\partial_i{\Bbb{G}}_{jk}+{\Bbb{G}}_{ki}\partial_j{\Bbb{K}}^i+{\Bbb{G}}_{ji}\partial_k{\Bbb{K}}^i=0.
\end{equation}
We will find the Killing vectors by solving these equations in
a particular order. Our approach is to first integrate out the $\eta$ and
$\bar{\eta}$ dependence.  Then the coefficients of the powers of
$\eta$ and $\bar{\eta}$ in the remaining equations must each be
zero, which determines the $\xi$ and $\bar{\xi}$ dependence. 

To start then, consider the simplest equations, which are ($j,k$) equal to ($\eta,\eta$),
($\eta,\bar{\eta}$) and ($\xi,\eta$). These say 
\begin{equation}\label{e:kill1}
\partial_\eta {\Bbb{K}}^{\bar{\xi}}=0
\end{equation}
\begin{equation}\label{e:kill2}
\partial_\eta {\Bbb{K}}^{\xi}-\partial_{\bar{\eta}} {\Bbb{K}}^{\bar{\xi}}=0
\end{equation}
\begin{equation}\label{e:kill3}
\partial_\xi {\Bbb{K}}^{\bar{\xi}}-\partial_{\eta} {\Bbb{K}}^{\bar{\eta}}=0.
\end{equation}
Differentiating equation (\ref{e:kill2}) with respect to $\eta$ and
noting equation (\ref{e:kill1}) we get $\partial_\eta\partial_\eta
{\Bbb{K}}^{\xi}=0$, with solution
\begin{equation}\label{e:kxi1}
{\Bbb{K}}^\xi=a_0(\xi,\bar{\xi},\bar{\eta})
    +a_1(\xi,\bar{\xi},\bar{\eta})\eta,
\end{equation}
where $a_0$ is complex-valued while $a_1$ is real-valued. Substituting
(\ref{e:kxi1}) into (\ref{e:kill1}) we get
$\partial_\eta\bar{a}_0+\bar{\eta}\partial_\eta a_1=0$. Looking at the
$\bar{\eta}$ dependence, both terms must be zero. That is,
$a_0=a_0(\xi,\bar{\xi})$ and $a_1=a_1(\xi,\bar{\xi})$ so that
\begin{equation}\label{e:kxi2}
{\Bbb{K}}^\xi=a_0(\xi,\bar{\xi}) +a_1(\xi,\bar{\xi})\eta.
\end{equation}
Substituting this into (\ref{e:kill3}) we get the equation $\partial_\eta
{\Bbb{K}}^{\bar{\eta}}=\partial_\xi\bar{a}_0+\bar{\eta}\partial_\xi a_1$,
which we can solve to get
\begin{equation}\label{e:keta1}
{\Bbb{K}}^\eta=b_0(\xi,\bar{\xi},\eta)+\left(\partial_{\bar{\xi}} a_0 
       +  \eta\partial_{\bar{\xi}} a_1\right)\bar{\eta}.
\end{equation}

At this stage (\ref{e:kxi2}) and (\ref{e:keta1}) solve 5 of the 10
Killing equations. The remaining equations are ($j,k$) equal to
($\xi,\bar{\eta}$), ($\xi,\xi$) and
($\xi,\bar{\xi}$). The first of these reads
\[
-(1+\xi\bar{\xi})(\partial_{\bar{\eta}}\bar{b}_0
  +\partial_\xi a_0 +  2\eta\partial_\xi a_1)
  +2a_0\bar{\xi}+2\bar{a}_0\xi+4a_1\xi\bar{\eta}=0.
\]
Now, the $\eta$ term tells us that $\partial_\xi a_1=0$ while its
conjugate implies that $a_1$ is constant. Integrating the remaining
equation with respect to $\bar{\eta}$ we find that
\[
\bar{b}_0=\bar{b}_1(\xi,\bar{\xi})
   -\left(\partial_\xi a_0-\frac{2(a_0\bar{\xi}
        +\bar{a}_0\xi)}{1+\xi\bar{\xi}}\right)\bar{\eta}
    +\frac{2a_1\xi}{1+\xi\bar{\xi}}\;\bar{\eta}^2.
\]

The $\eta$ and $\bar{\eta}$ dependence integrated out of the Killing
equations, we turn 
now to the $\xi$ and $\bar{\xi}$ dependence. Consider the Killing
equation with ($j,k$) equal to 
($\xi,\xi$). The coefficient of $\bar{\eta}^2$ is simply
$-2a_1$ which must therefore vanish. On the other hand the coefficient
of the term independent of both $\eta$ and $\bar{\eta}$ is
$\partial_\xi \bar{b}_1$ which must also be zero. Turning to the
$\eta$ term we get
\[
\partial_\xi\partial_\xi a_0
  +\frac{2\bar{\xi}}{1+\xi\bar{\xi}}\;\partial_\xi a_0=0,
\]
with solution
\[
a_0=c_0(\xi)+\frac{c_1(\xi)}{\xi(1+\xi\bar{\xi})}.
\]

We now turn to the final equation, ($j,k$) equal to
($\xi,\bar{\xi}$). The $\eta$ term of this gives
$-2(1+\xi\bar{\xi})\partial_{\bar{\xi}}\bar{c}_1+8\xi\bar{c}_1=0$,
which we solve to $c_1=c_2(1+\xi\bar{\xi})^4$, for some constant
$c_2$. The remaining part of the ($\xi,\bar{\xi}$) equation is
\[
\partial_\xi b_1-\frac{2\bar{\xi}}{1+\xi\bar{\xi}}\;b_1
  =\partial_{\bar{\xi}} \bar{b}_1
     -\frac{2\xi}{1+\xi\bar{\xi}}\;\bar{b}_1.
\]
The general solution to this is $b_1=c_3+c_4\xi-\bar{c}_3\xi^2$,
where $c_3$ is a complex constant while $c_4$ is a real constant. 

Returning finally to the ($\xi,\xi$) Killing equation, the $\eta$ term
tells us that $c_2=0$. The last Killing equation is then
\begin{equation}\label{e:final}
(1+\xi\bar{\xi})^2\partial_\xi\partial_\xi c_0
   -2\bar{\xi}(1+\xi\bar{\xi})\partial_\xi c_0
     +2\bar{\xi}\;^2c_0-2\bar{c}_0=0.
\end{equation}
Differentiating this three times with respect to $\bar{\xi}$ we
find that
$\partial_{\bar{\xi}}\partial_{\bar{\xi}}\partial_{\bar{\xi}}\bar{c}_0=0$,
so $c_0=d_0+d_1\xi+d_2\xi^2$. Substituting this back into
(\ref{e:final}) we get $d_0=\bar{d}_2$ and $d_1=-\bar{d}_1$.

A relabelling of constants gives the result.
\end{pf}

\section{Geodesics in ${\Bbb{T}}$}

Any curve in ${\Bbb{T}}$ gives a 1-parameter family of oriented lines
in ${\Bbb{R}}^3$. Classically, these are referred to as {\it ruled
  surfaces}, and we now determine the ruled surfaces that
correspond to the geodesics of ${\Bbb{G}}$.

\begin{Thm}  
A geodesic of ${\Bbb{G}}$ is either a plane or a helicoid in
${\Bbb{R}}^3$, the former in the case when the geodesic is null, the
latter when it is non-null. 
\end{Thm}
\begin{pf}
Let $c:[0,1]\rightarrow {\Bbb{T}}$ be a curve with tangent vector
\[
{\Bbb{X}}=\dot{\xi}\frac{\partial}{\partial \xi}
    +\dot{\eta}\frac{\partial}{\partial \eta}
   +\dot{\bar{\xi}}\frac{\partial}{\partial \bar{\xi}}
   +\dot{\bar{\eta}}\frac{\partial}{\partial \bar{\eta}}.
\]
By a translation and rotation we can set $\xi(0)=\eta(0)=0$, so that
the initial line $c(0)$ is the $x^3$-axis. We still
retain the freedom to rotate about, and translate along, the $x^3-axis$.

The geodesic equations (with arc-length or affine parameter $s$)
are

\begin{equation}\label{e:geo1}
\ddot{\xi}-\frac{2\bar{\xi}}{1+\xi\bar{\xi}}\dot{\xi}^2=0
\end{equation}
\begin{equation}\label{e:geo2}
\ddot{\eta}-\frac{4\bar{\xi}}{1+\xi\bar{\xi}}\dot{\xi}\dot{\eta}+\frac{2(\bar{\eta}+\bar{\xi}^2\eta)}{(1+\xi\bar{\xi})^2}\dot{\xi}^2=0,
\end{equation}
where a dot represents differentiation with respect to $s$, and we have
made use of the connection coefficients associated with ${\Bbb{G}}$. 

These have first integral
\begin{equation}\label{e:geo3}
\frac{2i}{(1+\xi\bar{\xi})^2}\left(
  \dot{\eta} \dot{\bar{\xi}}-\dot{\bar{\eta}} \dot{\xi}
   +\frac{2(\xi\bar{\eta}-\bar{\xi}\eta)}{1+\xi\bar{\xi}}\dot{\xi} \dot{\bar{\xi}}\right)=C_1,
\end{equation}
where $C_1\in{\Bbb{R}}$ vanishes iff the geodesic is null.
 
Suppose that $\dot{\xi}$ vanishes at $s=0$, then by (\ref{e:geo1}),
it vanishes along the whole geodesic. In addition, (\ref{e:geo2}) says
that $\eta$ is a linear function of the affine parameter and the ruled
surface is a plane. By (\ref{e:geo3}) the geodesic is null.

Suppose now that $\dot{\xi}(0)\neq 0$. Equation (\ref{e:geo1}) is the
geodesic equation for the round metric on $S^2$ and so the geodesic
projects to a great circle on $S^2$.  We can integrate the geodesic equations on the
sphere (with initial condition $\xi(0)=0$) to determine the evolution of $\xi$:
\begin{equation}\label{e:soln1}
\xi=\tan(C_2s)e^{i\theta},
\end{equation}
for constants $C_2\in{\Bbb{R}}$ and $\theta\in[0,2\pi)$. Substituting this in
(\ref{e:geo3}) we find that
\[
\eta e^{-i\theta}-\bar{\eta} e^{i\theta}
      =\frac{C_1s+C_3}{2iC_2\cos^2(C_2s)},
\]
for some real constant $C_3$. Equation (\ref{e:geo2}) now simplifies
to that of the forced harmonic oscillator:
\[
\frac{d^2}{ds^2}(\cos^2(C_2s)\eta)
   +4C_2^2\cos^2(C_2s)\eta =-C_2(C_1s+C_3)ie^{i\theta},
\] 
with solution
\[
\eta=\frac{C_4\cos(2C_2s)+C_5\sin(2C_2s)-(C_1s+C_3)i}{4C_2\cos^2(C_2s)}
    e^{i\theta}.
\]
Finally, since $\eta(0)=0$ we have $C_3=C_4=0$, so that
\begin{equation}\label{e:soln2}
\eta=\frac{C_5\sin(2C_2s)-C_1si}{4C_2\cos^2(C_2s)}
    e^{i\theta}.
\end{equation}
Equations (\ref{e:soln1}) and (\ref{e:soln2}) are
the general solution to the geodesic equations in ${\Bbb{T}}$ when
the initial line in ${\Bbb{R}}^3$ is the $x^3$-axis. The four real
constants remaining, namely $C_1,C_2,C_5$ and $\theta$, determine the
initial direction of the geodesic in ${\Bbb{T}}$. 

The associated ruled surface in ${\Bbb{R}}^3$ is
a helicoid when $C_1\neq 0$ and a plane when $C_1=0$. To see this we
can put it in standard position as follows. By a rotation about the
$x^3$-axis we can fix $\theta=0$, while a translation along the $x^3$-axis 
allows us to put $C_5=0$ ({\it cf}. equation (\ref{e:trans}) with
$a_1=c_1=0$ and $b_1=-C_5/2C_2$). The ruled surface can be explicitly
determined using (\ref{e:coord1}) and (\ref{e:coord2}) and the result
is
\[
x^1=t\sin(2C_2s) \qquad x^2=-\frac{C_1s}{2C_2} \qquad x^3=t\cos(2C_2s). 
\]
This is a helicoid for $C_1\neq 0$ and a plane for $C_1=0$, as claimed.
\end{pf}

Given two oriented lines $\gamma_1$ and $\gamma_2$, which neither
intersect nor are parallel, the geodesics in ${\Bbb{T}}$ joining them
consist of the helicoids in ${\Bbb{R}}^3$ containing them. The
${\Bbb{G}}$-distance between the lines is 
$-l^2/d$, where $l$ is the perpendicular distance between them in
${\Bbb{R}}^3$, and $d$ is the distance on ${\Bbb{P}}^1$ traversed by
the direction vectors of the ruling of the helicoid. This is maximised
(as appropriate for indefinite metrics) by the helicoid joining
$\gamma_1$ and $\gamma_2$ with the minimum number of turns. The
multiplicity of a non-maximising geodesic in ${\Bbb{G}}$ is given by
the number of turns of the helicoid.

\section{Surfaces in ${\Bbb{T}}$}

Consider now a {\it line congruence} $\Sigma$, that is, a 2-parameter
family of oriented lines in ${\Bbb{R}}^3$. Equivalently, 
this is a mapping $f:\Sigma \rightarrow {\Bbb{T}}$ of a surface
$\Sigma$ into ${\Bbb{T}}$. 

Away from crossings of the lines we can adapt a local null
frame to $\Sigma$ by aligning $e_{(0)}$ with the direction of the
lines. The complex {\it spin coefficients} of the congruence are defined by
\[
\Gamma_{mnp}=e_{(n)}^ie_{(p)}^j\nabla_je_{(m)i},
\]
where $\nabla$ is the Euclidean covariant derivative and the indices
$m$, $n$, $p$ range over $0$, $+$, $-$. Breaking covariance, introduce
the complex optical scalars 
\[
\Gamma_{+0-}=\rho \qquad\qquad \Gamma_{+0+}=\sigma. 
\]

The complex scalar functions $\rho$ and $\sigma$ describe
the first order geometric behaviour of the congruence of lines. 
 The real part of $\rho$ is the {\it
  divergence}, the imaginary part $\lambda$ is the {\it twist} and
$\sigma$ is the {\it shear} of the congruence (see \cite{guil2} and
\cite{par} for further details). 

\begin{Prop}
A surface $\Sigma$ in ${\Bbb{T}}$ is Lagrangian with respect to the
symplectic structure $\Omega$ if and only if the
associated congruence is integrable (twist-free) i.e. there exists a
surface $S$ in ${\Bbb{R}}^3$ orthogonal to the line congruence.
\end{Prop}
\begin{pf}
A surface $\Sigma$ is Lagrangian iff $\Omega$ pulled back to $\Sigma$
vanishes. Suppose that $\Sigma$ is given parametrically by
$f:\Sigma\rightarrow {\Bbb{T}}$, $\nu\mapsto
(\xi(\nu,\bar{\nu}),\eta(\nu,\bar{\nu}))$. Then pulling $\Omega$
back to $\Sigma$, using (\ref{e:riemsymp}):
\[
f_*\Omega=4{\Bbb{R}}\mbox{e}\left(\left(\partial\eta\bar{\partial}\bar{\xi}
               +\partial\bar{\eta}\bar{\partial}\xi
          -2\frac{\bar{\xi}\eta}{1+\xi\bar{\xi}}(\partial\xi\bar{\partial}\bar{\xi}-\bar{\partial}\xi\partial\bar{\xi})\right)\frac{d\nu \wedge d\bar{\nu}}{(1+\xi\bar{\xi})^2}\right),
\]
where $\partial=\frac{\partial}{\partial \nu}$.

This is precisely the expression for the twist of a line congruence
derived in \cite{gak2}, which 
vanishes iff there exists surfaces $S$ in ${\Bbb{R}}^3$ orthogonal to
the line congruence (by Frobenius' theorem). 
\end{pf}

The canonical 1-form $\Theta$ pulled back
to a Lagrangian surface $\Sigma$ is closed. Thus it defines an element in the
real cohomology $H^1(\Sigma,{\Bbb{R}})$. Locally, $\Theta=dr$, where
$r$ is a real function on $\Sigma$. In fact, the surfaces in
${\Bbb{R}}^3$ orthogonal to the line congruence are obtained locally
by substituting $r=r(\nu,\bar{\nu})$ in (\ref{e:coord1}) and
(\ref{e:coord2}). The evaluation of $\Theta$ on a closed curve in
$\Sigma$ is equal to the jump of $r$ as one goes around a curve on the
orthogonal surfaces in
${\Bbb{R}}^3$.

\vspace{0.1in}
\noindent{\bf Example 1.}
Consider the Lagrangian torus in ${\Bbb{T}}$ given by 
\[
\xi=\tan\phi \;e^{i\theta} \qquad \qquad
\eta=\pm a(1-b\tan^2\phi)e^{i\theta},
\]
for $\theta,\phi\in[0,\pi)$. For $b=1$ this is the normal congruence
to the rotationally symmetric torus in ${\Bbb{R}}^3$ with core radius $2a$.
For $b\neq 1$ it is a rotationally symmetric torus ``torn'' along an
equilateral. The jump in $r$ as one traverses a meridian is equal to
$2\pi a(1-b)$. 

\vspace{0.2in}

Turning to the metric ${\Bbb{G}}$ we have:

\begin{Thm}
The induced metric on a surface $\Sigma$ in ${\Bbb{T}}$  is 
Lorentzian (totally null, Riemannian) if and only if
$\lambda^2-|\sigma|^2<0 \;(=0,>0)$,
where $\lambda$ and $\sigma$ are the twist and shear of the line
congruence $\Sigma$.
\end{Thm}
\begin{pf}
Again, suppose that $\Sigma$ is given parametrically by
$f:\Sigma\rightarrow {\Bbb{T}}$, $\nu\mapsto
(\xi(\nu,\bar{\nu}),\eta(\nu,\bar{\nu}))$. Then pulling back
${\Bbb{G}}$ to $\Sigma$:
\begin{align}
f_*{\Bbb{G}}=&{\Bbb{I}}\mbox{m}\left(\frac{4}{(1+\xi\bar{\xi})^2}\left[\left( \partial\eta\partial\bar{\xi}
       -\partial\bar{\eta}\partial\xi
         +2\frac{(\xi\bar{\eta}-\bar{\xi}\eta)}{1+\xi\bar{\xi}}\partial\xi\partial\bar{\xi}\right)d\nu^2\right.\right.\nonumber\\
&\qquad\qquad+\left.\left.\left(\partial\eta\bar{\partial}\bar{\xi}
       +\bar{\partial}\eta\partial\bar{\xi}
   -2\frac{\bar{\xi}\eta}{1+\xi\bar{\xi}}(\partial\xi\bar{\partial}\bar{\xi}-\bar{\partial}\xi\partial\bar{\xi})\right)d\nu d\bar{\nu}\right]\right).\nonumber
\end{align}
By the expressions derived in \cite{gak2} this is equivalent to
\begin{align}\label{e:indmet}
f_*{\Bbb{G}}=&{\Bbb{I}}\mbox{m}\left(\frac{4}{K(1+\xi\bar{\xi})^2}\left[
   \left((\rho-\bar{\rho})\partial\xi\partial\bar{\xi}+\sigma(\partial\xi)^2-\bar{\sigma}(\partial\bar{\xi})^2\right)d\nu^2\right.\right.\nonumber\\
&\qquad\qquad+\left.\left.\left(\rho(\partial\xi\bar{\partial}\bar{\xi}+\bar{\partial}\xi\partial\bar{\xi})+2\sigma\partial\xi\bar{\partial}\xi\right)d\nu d\bar{\nu}\right]\frac{}{}\right),
\end{align}
where $K$ is the (generalised) curvature of the congruence. The
determinant of this matrix is proportional to $\lambda^2-|\sigma|^2$,
where $\lambda$ is the imaginary part of $\rho$, and the theorem follows.
\end{pf}

\vspace{0.1in}
\noindent{\bf Example 2.}
Suppose $\Sigma$ is a Lagrangian surface so that
$\lambda$ is zero, the above theorem implies that the induced metric is
either Lorentz or totally null. In this case the geometric
significance of the metric is as follows. 

Equation (\ref{e:indmet}) implies
that the null directions of the metric are given by the argument of
the shear. According to the results in \cite{gak2} these are the
principal directions of the orthogonal surface $S$ in
${\Bbb{R}}^3$. Thus the null curves on $\Sigma$ correspond to the lines of
curvature on $S$. Similarly, the totally null points on $\Sigma$ are the
shear-free lines in the congruence and these are the umbilical points
on $S$.

\vspace{0.1in}

\noindent{\bf Example 3.}
Suppose $\Sigma$ is holomorphic. Then the induced metric is either
positive definite or totally null. To construct examples, note that a
vector field on $S^2$ gives rise to a surface $\Sigma$ in ${\Bbb{T}}$
and consider the vector field generated on the unit
sphere in ${\Bbb{R}}^3$ by a rotation  about
the $x^3$-axis. This is the line congruence given by the global
holomorphic section $\xi\mapsto (\xi,\eta=-bi\xi)$ for
$b\in{\Bbb{R}}^3$. The congruence is twisting everywhere except along the
equator $|\xi|=1$. The induced metric is 
\[
{\Bbb{G}}_\Sigma
    =\frac{4b(1-\xi\bar{\xi})}{(1+\xi\bar{\xi})^3}d\xi d\bar{\xi},
\]
which is positive definite except along the
equator, where it is totally null. The Northern and Southern
hemispheres are overtwisted discs \cite{eliah}. 

\vspace{0.1in}

\noindent{\bf Example 4.}
The only line congruences that are both Lagrangian and
holomorphic are the normal congruences to round spheres and flat
planes in ${\Bbb{R}}^3$. 

 \vspace{0.1in}

\section{The Keller-Maslov Index}

Given a surface $\Sigma$ in a symplectic 4-manifold
with compatible tamed complex structure, it is well-known that one can
define the Keller-Maslov index of a totally real curve on $\Sigma$
\cite{abk}.  We will now mirror this construction on Lagrangian
surfaces in ${\Bbb{T}}$, where we will use the complex structure
${\Bbb{J}}$, despite the fact that it is not tamed by the symplectic
structure $\Omega$.

Let $\Lambda_\gamma$ be the Grassmanian of 2-planes at
$\gamma\in{\Bbb{T}} $. The Gauss map of $\Sigma$ takes a point $\gamma\in\Sigma$
to an element $T_\gamma\Sigma\in\Lambda_\gamma$. Let the subgroup of
matrices in $\mbox{GL}(2,{\Bbb{C}})$ which preserve $T_\gamma\Sigma$
be denoted by $G_\gamma$. 

A point $\gamma\in\Sigma$ is {\it complex} if
${\Bbb{J}}$ preserves $T_\gamma\Sigma$ and such points, when isolated,
have an associated multiplicity. In the case of ${\Bbb{T}}$, a point
$\gamma$ is complex with respect to ${\Bbb{J}}$ if and only if the
associated line in the congruence has vanishing shear \cite{gak2}. 

Using the indefinite metric ${\Bbb{G}}$ we reduce the group
$\mbox{GL}(2,{\Bbb{C}})$ to $U(1,1)=\mbox{PSL}(2,{\Bbb{R}})\times\mbox{S}^1$ and, as proved in Theorem 3, at a
Lagrangian point, the metric is either Lorentz or totally null, the
latter occurring only at complex points. Thus, away from complex points,
$G_\gamma$ can be reduced to $O(1,1)$ and $\Lambda_\gamma/G_\gamma$ is
isomorphic to $S^1$. The Keller-Maslov index of a closed oriented
totally real curve on a Lagrangian surface 
$\Sigma$ is the winding number of $T_\gamma\Sigma$  within
$\Lambda_\gamma$ . 

\begin{Thm}
Suppose that the complex points on a Lagrangian surface
$\Sigma\subset{\Bbb{T}}$ are isolated. Then the Keller-Maslov index of
a closed oriented curve on $\Sigma$ is
equal to the number of shear-free lines inside the curve (counted with
multiplicity). 
\end{Thm}

\begin{pf}

A complex point is a shear-free line, or, as
$\Sigma$ is Lagrangian, an
umbilical point on the orthogonal surface $S$ in ${\Bbb{R}}^3$
\cite{gak2}.  By the definition above, the 
Keller-Maslov index measures the rotation of the null directions on
the Lorentz surface $\Sigma$, or equivalently, the rotation of the
principal foliation on $S$. 

Thus the index counts (with multiplicity) the number of shear-free 
lines inside the curve on $\Sigma$.
\end{pf}

\end{document}